\documentclass[11pt,reqno]{amsart}
\usepackage{graphicx}
\usepackage{amsfonts}
\usepackage{mathrsfs}
 \usepackage{times}
 \usepackage{mathptmx}
\usepackage{amsmath}
\usepackage{txfonts}
\usepackage{amssymb}
\usepackage{amsmath,amsthm}
\usepackage{CJK,amsfonts}
\usepackage[body={7.5in,8.75in}, top=1.2in, right=0.9in,left=0.9in]{geometry}
\usepackage[colorlinks=true, linkcolor=red, citecolor=blue]{hyperref}
\allowdisplaybreaks[4] \numberwithin{equation}{section}

\newtheorem{thm}{Theorem}

\newtheorem{corollary}[thm]{Corollary}
\newtheorem{lemma}[thm]{Lemma}
\newtheorem{remark}[thm]{Remark}
\newtheorem{proposition}[thm]{Proposition}

\def\squarebox#1{\hbox to #1{\hfill\vbox to #1{\vfill}}}

\newcommand{\qbinomial}[3]{\mbox{$
\biggl[ 
\begin{array}{c}
#1\\
 #2
\end{array}\biggr]_{
\!{#3}}$}}

\newcommand{\be}{\begin{equation}}
\newcommand{\ee}{\end{equation}}
\newcommand{\bea}{\begin{eqnarray}}
\newcommand{\eea}{\end{eqnarray}}
\newcommand{\bd}{\begin{displaymath}}
\newcommand{\ed}{\end{displaymath}}
\begin{document}
\begin{CJK*}{GBK}{song}
\title[Certain  generating functions  for Cigler's polynomials]{Certain  generating functions  for Cigler's polynomials}
\author{  Sama Arjika${}^{1,2}$ }
\dedicatory{\textsc}
\thanks{${}^1$Department of Mathematics and Informatics, University of Agadez, Niger.  ${}^2$International Chair of Mathematical Physics and Applications (ICMPA-UNESCO Chair),
University of Abomey-Calavi, Post Box 072, Cotonou 50,  Benin}
\thanks{Email:  rjksama2008@gmail.com.}



\begin{abstract}
In this   paper, we   use the  homogeneous $q$-operators [J.  Difference  Equ. Appl. {\bf20 } (2014),  837--851.] to derive    Rogers formulas,    extended Rogers formulas and Srivastava-Agarwal   type bilinear generating functions  for Cigler's polynomials [J. Difference Equ. Appl. {\bf 24} (2018),   479--502.].   Finally,  we also derive two interesting   transformation formulas between  ${}_2\Phi_1, \,  {}_2\Phi_2$  and ${}_3\Phi_2$.
\end{abstract}

\maketitle

\section{\bf Introduction} 
Wang and Cao \cite{CAO2018} presented two extensions of Cigler's
polynomials together with their several generating functions for these extended
Cigler's polynomials by using the method of homogeneous $q$-difference equations.
In this paper, by employing the homogeneous $q$-operators, we aim to establish
more generalized generating functions for the extended Cigler's polynomials than
those in Wang and Cao \cite{CAO2018}.

The  $q$-Laguerre polynomials are defined by \cite{Koekock} 
\be 
\label{def00}
\mathcal{L}_n^{(\alpha)}(x;q)=\frac{(q^{\alpha+1};q)_n}{(q;q)_n}
\sum_{k=0}^n  \frac{q^{({}^{k}_{2})}(q^{-n};q)_k}{(q^{\alpha+1},q;q)_{k}}  (xq^{n+\alpha+1})^k 
\ee 
which  belong to  the Askey-scheme of basic hypergeometric orthogonal polynomials according to Koekoek and Swarttouw \cite[Eq. (3.21.1)]{Koekock}. They appear in several branches of  mathematics   and physics, and their generalization arise in many  applications (see for details, \cite{Andrews74,Atakishiyeva97,Atakishiyev03,Chung98,Coulembier10,Micu05}), 
 such as (for example) quantum group,   $q$-harmonic oscillator and coding theory, and so on.

Cigler studied  $q$-Laguerre polynomials  \cite[Eq. (30)]{Cigler81}
\be 
\label{def0}
l_n^{(\alpha)}(x)=\frac{1}{q^{n^2+\alpha n}}
\sum_{k=0}^n {n+\alpha\,\atopwithdelims []\,n-k\,}_{q} \frac{(q;q)_n}{(q;q)_{k}}(-1)^kq^{({}^{n-k}_{\,\,\,2})} x^{n-k}
\ee 
and derived  important results.

In 2013, Cao and Niu \cite{CAO2016} introduced two  extensions of  Cigler's  polynomials 
\be 
\label{def0}
\mathcal{C}_n^{(\alpha)}(x,b)= 
\sum_{k=0}^n (-1)^kq^{({}^k_2)} {n+\alpha\,\atopwithdelims []\,k\,}_{q} \frac{(q;q)_n}{(q;q)_{k}}  b^k x^{n-k}
\ee 
and
\be 
\label{def01}
\mathcal{D}_n^{(\alpha)}(x,b)= 
\sum_{k=0}^n q^{k^2-nk} {n+\alpha\,\atopwithdelims []\, k\,}_{q} \frac{(q;q)_n}{(q;q)_{k}}b^k x^{n-k}
\ee 
 as solution of $q$-difference equations and deduced their generating functions. 

 Recently, Wang and Cao \cite{CAO2018}  introduced two extensions of Cigler's polynomials
 \be 
\label{def02}
\mathcal{C}_n^{(\alpha-n)}(x,y,b)= 
\sum_{k=0}^n (-1)^kq^{({}^k_2)} { \alpha\,\atopwithdelims []\,k\,}_{q}b^k \frac{(q;q)_n}{(q;q)_{n-k}}   p_{n-k}(x,y) 
\ee 
and
\be 
\label{def03}
\mathcal{D}_n^{(\alpha-n)}(x,y,b)= 
\sum_{k=0}^nq^{({}^k_2)}  { \alpha\,\atopwithdelims []\, k\,}_{q}b^k\frac{(q;q)_n}{(q;q)_{n-k}}\left[(-1)^{n+k} q^{-({}^n_2)}  p_{n-k}(y,x)\right],
\ee 
and derived  the following important results by the homogeneous $q$-difference equations. 
 \begin{proposition}
  \cite[Theorem 5]{CAO2018}
\label{dgend} For  $\alpha\in\mathbb{R}$ and $s\in\mathbb{N}$, we have:
 \be 
\label{exdddgen}
\sum_{n=0}^\infty  \mathcal{C}_{n+s}^{(\alpha-n-s)}(x,y,b)\frac{  t^{n+s}}{(q;q)_n} 
= (bt;q)_\alpha\frac{ (yt;q)_\infty}{ (xt;q)_\infty} {}_{3}\Phi_2\left[\begin{array}{r}q^{-s},xt, q^{ \alpha}bt;
 \\\\
yt, bt;
 \end{array} 
q; q\right], \,\, |xt|<1,
 \ee
 \be 
\label{sexdddgen}
\sum_{n=0}^\infty \mathcal{D}_{n+s}^{(\alpha-n-s)}(x,y,b)(-1)^{n+s}q^{ ({}^{n+s}_{\,\,\,2})} \frac{  t^{n+s}}{(q;q)_n}  
=(bt;q)_\alpha\frac{ (xt;q)_\infty}{ (yt;q)_\infty}  {}_{3}\Phi_2\left[\begin{array}{r}q^{-s},yt, q^{ \alpha}bt;
 \\\\
xt, bt;
 \end{array} 
q; q\right], \,\, |yt|<1. 
 \ee
 \end{proposition}
 \begin{proposition}
  \cite[Corollary 6]{CAO2018}
\label{Corl}
For $\alpha\in\mathbb{R}$,  we have:  
\be  
\label{ler}
\sum_{n=0}^\infty  \mathcal{C}_{n }^{(\alpha-n)}(x,y,b)\frac{  t^{n }}{(q;q)_n} 
= \frac{ (yt, bt;q)_\infty}{ (xt,q^{ \alpha}bt;q)_\infty},\quad \max\{|xt|,|q^\alpha bt|\}<1,\ee
 \begin{equation}
\label{lerr}
\sum_{n=0}^\infty(-1)^{n}q^{ ({}^{n}_{2})}\mathcal{D}_{n}^{(\alpha-n)}(x,y,b) \frac{  t^{n }}{(q;q)_n} 
=\frac{ (xt, bt;q)_\infty}{ (yt,q^{ \alpha}bt;q)_\infty},\quad  \max\{|yt|,|q^\alpha bt|\}<1.
 \end{equation}
 \end{proposition}
 \begin{proposition}
  \cite[Theorem 11]{CAO2018}
\label{erTA1} For $\alpha\in\mathbb{R}$,  we have:
 \be 
\label{gs}
  \sum_{n=0}^\infty  \mathcal{C}_{n}^{(\alpha-n )}(x,y,b)  (\lambda;q)_n \frac{  t^n }{(q;q)_n}  
= \frac{(\lambda,yt, bt;q)_\infty}{(xt, q^{ \alpha}bt;q)_\infty}  {}_{3}\Phi_2\left[\begin{array}{r}xt, q^{ \alpha}bt,0; 
 \\\\
yt, bt;
 \end{array} 
q; \lambda\right],
 \ee
where $ \max\{ |xt|,|\lambda|,|q^\alpha bt|\}<1;$
 \be 
\label{cc1sums}
 \sum_{n=0}^\infty(-1)^nq^{({}^n_2)}\,   \mathcal{D}_{n}^{(\alpha-n )}(x,y,b)  (\lambda;q)_n \frac{ t^n }{(q;q)_n}  
\\= \frac{(\lambda,xt, bt;q)_\infty}{(yt, q^{ \alpha}bt;q)_\infty}  {}_{3}\Phi_2\left[\begin{array}{r}yt, q^{ \alpha}bt,0;
 \\\\
xt, bt;
 \end{array}
q; \lambda\right],
 \ee
where $ \max\{ |yt|,|\lambda|,|q^\alpha bt|\}<1.$
 \end{proposition}
 \begin{remark}
  For  $s=0$,  the assertions  (\ref{exdddgen}) and (\ref{sexdddgen}) of Proposition \ref{dgend}    reduce to the assertions (\ref{ler}) and (\ref{lerr}).  
  \end{remark}
  \begin{remark}
   For  $\lambda=0$, the assertions  (\ref{gs}) and (\ref{cc1sums})   of  Proposition \ref{erTA1} reduce to the assertions (\ref{ler}) and (\ref{lerr}).  
  \end{remark}
  
In this paper,   motivated by Wang and Cao's  results \cite{CAO2018},  we aim to establish
more generalized generating functions for the extended Cigler's polynomials $\mathcal{C}_{n }^{(\alpha-n)}(x,y,b)$ and $\mathcal{D}_{n }^{(\alpha-n)}(x,y,b)$.

 Our main results   are stated as:
\begin{thm}
 {\rm $\big[$Rogers type formulas for $\mathcal{C}_{n }^{(\alpha-n)}(x,y,b)$ and $\mathcal{D}_{n }^{(\alpha-n)}(x,y,b)\big]$} 
\label{exgend} 
For $\alpha \in\mathbb{R},$ the following Rogers type formulas  hold  true for  Cigler's   polynomials $\mathcal{C}_{n }^{(\alpha-n)}(x,y,b)$ and  $\mathcal{D}_{n }^{(\alpha-n)}(x,y,b)$$:$
 \be 
\label{exgen}
 \sum_{n=0}^\infty\sum_{m=0}^\infty   \mathcal{C}_{n+m}^{(\alpha-n-m)}(x,y,b)\frac{  t^{n }}{(q;q)_n}\frac{  s^{m }}{(q;q)_m} 
=\frac{ (ys,bs;q)_\infty}{ (t/s,xs,q^{\alpha} bs;q)_\infty}    {}_{4}\Phi_3\left[\begin{array}{r}xs,q^{\alpha} bs,0,0;
 \\\\
qs/t,ys, bs;
 \end{array} 
q; q\right], 
 \ee
where $\max\{|t/s|,|xs|,|q^\alpha  bs|\}<1$;
\be 
\label{2exgen}
 \sum_{n=0}^\infty\sum_{m=0}^\infty   (-1)^{n+m}q^{({}^{n+m}_{\,\,\,\,\,2})}\mathcal{D}_{n+m}^{(\alpha-n-m)}(x,y,b) \frac{  t^{n }}{(q;q)_n}\frac{  s^{m }}{(q;q)_m}
=\frac{ (xs,bs;q)_\infty}{ (t/s,ys,q^{\alpha} bs;q)_\infty}    {}_{4}\Phi_3\left[\begin{array}{r}ys,q^{\alpha} bs,0,0;
 \\\\
qs/t,xs, bs;
 \end{array} 
q; q\right],
\ee
where $ \max\{|t/s|,|ys|,|q^\alpha bs|\}<1$.
 \end{thm}
 \begin{thm}
  {\rm $\big[$Extended Rogers type formula for $\mathcal{C}_{n }^{(\alpha-n)}(x,y,b)$ and $\mathcal{D}_{n }^{(\alpha-n)}(x,y,b)\big]$} 
\label{Tsxt}
 For $\alpha \in\mathbb{R},$ the following extended Rogers type formulas   hold  true for  Cigler's   polynomials $\mathcal{C}_{n }^{(\alpha-n)}(x,y,b)$ and  $\mathcal{D}_{n }^{(\alpha-n)}(x,y,b)$$:$
 \begin{multline}
\label{gextend}
\sum_{n=0}^\infty \sum_{m=0}^\infty \sum_{k=0}^\infty  \mathcal{C}_{n+m+k}^{(\alpha-n-m-k)}(x,y,b)\frac{ t^n\,s^m\omega^k}{(q;q)_{n+m}(q;q)_{m} (q;q)_k}\\
=  \frac{ (y\omega,b\omega;q)_\infty}{ (s/t,t/\omega,x\omega,q^{ \alpha} b\omega;q)_\infty} 
     {}_{4}\Phi_3\left[\begin{array}{r}x\omega,q^{ \alpha} b\omega,0,0; 
 \\\\
y\omega, b\omega,q\omega/t;
 \end{array} 
q; q\right],
\end{multline}
where $ \max\{|s/t|,|t/\omega|,|x\omega|,|q^\alpha b\omega|\}<1; $
\begin{multline}
\sum_{n=0}^\infty \sum_{m=0}^\infty \sum_{k=0}^\infty \mathcal{D}_{n+m+k}^{(\alpha-n-m-k)}(x,y,b)\frac{ (-1)^{n+m+k}q^{({}^{n+m+k}_{\,\,\,\,\,\,\,\,\,2})} t^n\,s^m\omega^k}{(q;q)_{n+m}(q;q)_{m} (q;q)_k}\\
=  \frac{ (x\omega,b\omega;q)_\infty}{ (s/t,t/\omega,y\omega,q^{ \alpha} b\omega;q)_\infty} 
     {}_{4}\Phi_3\left[\begin{array}{r}y\omega,q^{ \alpha} b\omega,0,0; 
 \\\\
x\omega, b\omega,q\omega/t;
 \end{array} 
q; q\right],\label{2gextend}
\end{multline}
where $ \max\{|s/t|,|t/\omega|,|y\omega|,|q^\alpha b\omega|\}<1.$
 \end{thm}
  \begin{thm} {\rm $\big[$Srivastava-Agarwal type bilinear generating function for $\mathcal{C}_{n }^{(\alpha-n)}(x,y,b)$ and $\mathcal{D}_{n }^{(\alpha-n)}(x,y,b)\big]$} 
  \label{TA1}
For $\alpha,\,\beta\in\mathbb{R},$ we have:  
\be
\label{1sums}
\sum_{n=0}^\infty  \phi_n^{(\alpha)}(x|q) \mathcal{C}_{n}^{(\beta-n )}(u,v,b)  \frac{ \,t^n}{(q;q)_n}  
= \frac{ (\alpha x,v t,bt;q)_\infty}{ (x,u t,q^{\beta} bt;q)_\infty}    {}_{4}\Phi_3\left[\begin{array}{r}\alpha,u t,q^{\beta} bt,0;
 \\\\
q/x,v t, bt;
 \end{array} 
q; q\right],
\ee
where $ \max\{ |x|,|u t|,|q^{\beta} bt|\}<1;$
\begin{multline} 
\label{2sums}
\sum_{n=0}^\infty(-1)^nq^{({}^{n+1}_{\,\,\,\,\,2})}  \psi_n^{(\alpha)}(x|q) \mathcal{D}_{n}^{(\beta-n )}(u,v,b)  \frac{ t^n}{(q;q)_n}\\
\qquad\qquad= \frac{   (q/x,uxtq,  bxtq; q)_\infty }{ (\alpha q,v  xtq,  bxtq^{1+\beta};q)_\infty}  {}_{3}\Phi_3\left[\begin{array}{rr}1/(\alpha x),1/(uxt),  1/(bxt);\\
 \\
q/x,1/(v  xt),  q^{-\beta }/(bxt);
 \end{array} 
q; \frac{\alpha u q^{1-\beta }}{v} \right],
\end{multline}
where $ \max\{|\alpha q|,|v  xt|,  |bxtq^{\beta+1}|\}<1.$
 \end{thm} 
\begin{remark}
For $k=0$, Theorem \ref{Tsxt} reduces to Theorem \ref{exgend}.  
\end{remark}
         \begin{corollary} 
  \label{CTA1}
For $\alpha,\;\beta\in\mathbb{R},$ we have:  
\be 
\label{1s}
\sum_{n=0}^\infty  \phi_n^{(\alpha)}(x|q) \mathcal{C}_{n}^{(\beta-n )}(u,b)  \frac{ \,t^n}{(q;q)_n}  
\\= \frac{ (\alpha x,bt;q)_\infty}{ (x,u t,q^{\beta} bt;q)_\infty}    {}_{3}\Phi_2\left[\begin{array}{r}\alpha,u t,q^{\beta} bt;
 \\\\
q/x, bt;
 \end{array} 
q; q\right],\ee
where $\max\{ |x|,|u t|,|q^{\beta} bt|\}<1;$
\be 
\label{2s}
\sum_{n=0}^\infty(-1)^nq^{({}^{n+1}_{\,\,\,\,\,2})}  \psi_n^{(\alpha)}(x|q) \mathcal{D}_{n}^{(\beta-n )}(u,b)  \frac{ t^n}{(q;q)_n}\\
= \frac{   (q/x,uxtq,  bxtq; q)_\infty }{ (\alpha q,  bxtq^{1+\beta};q)_\infty}  {}_{3}\Phi_2\left[\begin{array}{rr}1/(\alpha x),1/(uxt),  1/(bxt);
 \\\\
q/x,  q^{-\beta }/(bxt);
 \end{array} 
q;  \alpha u xtq^{1-\beta } \right],
\ee
where $ \max\{|\alpha q|,  |bxtq^{\beta+1}|\}<1.$
 \end{corollary} 
    \begin{remark}
   For $v=0$, the assertions  (\ref{1sums}) and (\ref{2sums}) of  Theorem \ref{TA1}, reduce to the assertions   (\ref{1s}) and (\ref{2s}) .
    \end{remark}
    
 The rest of paper is organized as follows. In Section \ref{prel},  we present   notations and give some $q$-operator identities.  In Section \ref{prel3}, we use the homogeneous $q$-operators to derive   Rogers formulas and extended Rogers formulas.    In Section \ref{qdiffff}, we give  Srivastava-Agarwal type bilinear generating functions for   generalized Cigler's  polynomials.  As   an application of Srivastava-Agarwal type generating functions,    we deduce two   interesting   transformation formulas between  ${}_2\Phi_1, \,  {}_2\Phi_2$  and ${}_3\Phi_2$ in    Section \ref{q5}.  We end by the concluding remarks in   Section \ref{conclusion}.
 
 \section{Notations and lemmas}
\label{prel}
In this section, we adopt the common   notation and terminology for basic hypergeometric  series as in Refs.  \cite{GasparRahman,Koekock}.    Throughout this paper, we assume that  $q$ is a fixed nonzero real or complex number and $|q|< 1$.  The $q$-shifted factorial  and its compact  factorial  are defined \cite{GasparRahman,Koekock}, respectively by:
\begin{equation}
(a;q)_0:=1,\quad  (a;q)_{n} =\prod_{k=0}^{n-1} (1-aq^k),   \; (a;q)_{\infty}:=\prod_{k=0}^{\infty}(1-aq^{k})
\end{equation}
and 
 $ (a_1,a_2, \ldots, a_r;q)_m=(a_1;q)_m (a_2;q)_m\cdots(a_r;q)_m,\; m\in\{0, 1, 2\cdots\}$. \\
 We will use frequently  the following relation
 \begin{equation}
 \label{usu}
(aq^{-n};q)_n=(q/a;q)_n(-a)^nq^{-n-({}^n_2)}. 
\end{equation}
  The generalized $q$-binomial coefficient is defined as \cite{GasparRahman}
\be
\label{qb}
 {\,\alpha\,\atopwithdelims []\,k\,}_{q}=\frac{(q^{-\alpha};q)_k}{(q;q)_k}(-1)^kq^{\alpha k-({}^k_2)},\,\,\alpha \in\mathbb{C}.
\ee
\par 
 Here, in our present investigation, we are mainly concerned  
with the Cauchy polynomials $p_n(x,y)$ as given
below (see \cite{Chen2003} and \cite{GasparRahman}):
\bea
\label{def}
 p_n(x,y):=(x-y)( x- qy)\cdots ( x-q^{n-1}y) =(y/x;q)_n\,x^n
\eea
 which has the following  generating function \cite{Chen2003}
\be
\label{gener}
\sum_{n=0}^{\infty} p_n(x,y)
\frac{t^n }{(q;q)_n} = 
\frac{(yt;q)_\infty}{(xt;q)_\infty}.
\ee
 The generating   function (\ref{gener}) is also the   homogeneous version  of the Cauchy identity or the $q$-binomial theorem   given by \cite{GasparRahman}
\be
\label{putt}
\sum_{k=0}^{\infty} 
\frac{(a;q)_k }{(q;q)_k}z^{k}={}_{1}\Phi_0\left[\begin{array}{c}a;
 \\\\
-;
 \end{array} 
q;z\right]= 
\frac{(az;q)_\infty}{(z;q)_\infty},\quad |z|<1,  
\ee
where the    basic or $q$-hypergeometric function    in the variable $z$ (see  Slater \cite[Chap. 3]{SLATER},  Srivastava and Karlsson   \cite[p. 347, Eq. (272)]{SrivastaKarlsson}   for details) is defined as:
 $$
{}_{r}\Phi_s\left[\begin{array}{r}a_1, a_2,\ldots, a_r;
 \\\\
b_1,b_2,\ldots,b_s;
 \end{array} 
q;z\right]
 =\sum_{n=0}^\infty\Big[(-1)^n q^{({}^n_2)}\Big]^{1+s-r}\,\frac{(a_1, a_2,\ldots, a_r;q)_n}{(b_1,b_2,\ldots,b_s;q)_n}\frac{ z^n}{(q;q)_n},
$$
 when $r>s+1$. Note that, for $r=s+1$, we have:
$$
{}_{r+1}\Phi_r\left[\begin{array}{r}a_1, a_2,\ldots, a_{r+1};
 \\\\
b_1,b_2,\ldots,b_r;
 \end{array} q; z\right]
 =\sum_{n=0}^\infty \frac{(a_1, a_2,\ldots, a_{r+1};q)_n}{(b_1,b_2,\ldots,b_r;q)_n}\frac{ z^n}{(q;q)_n}.
$$
Putting    $a=0$, the relation (\ref{putt}) becomes  Euler's identity  \cite{GasparRahman}
\be
\label{q-expo-alpha}
  \sum_{k=0}^{\infty} \frac{ z^{k}}{(q;q)_k}=\frac{1}{(z;q)_\infty},\quad |z|<1
\ee
and its inverse relation  \cite{GasparRahman}
\be
\label{q-Expo-alpha}
 \sum_{k=0}^{\infty}  \frac{(-1)^kq^{ ({}^k_2)
}\,z^{k}}{(q;q)_k}=(z;q)_\infty.
\ee
 \par
We remark in passing that, in a recently-published 
survey-cum-expository review article, the so-called $(p,q)$-calculus 
was exposed to be a rather trivial and inconsequential variation of 
the classical $q$-calculus, the additional parameter $p$ being redundant 
or superfluous (see, for details, \cite[p. 340]{HMS-ISTT2020}).\par
  Chen {\it et al.} \cite{Chen2003} introduced  homogeneous $q$-difference operator   $D_{xy}$,  Saad and  Sukhi \cite{Saadsukhi} introduced another homogeneous $q$-difference operator ${\theta}_{xy}$  as
\be
\label{deffd}
  D_{xy}\big\{f(x,y)\}:=\frac{f(x,q^{-1}y)-f( qx, y)}{x-q^{-1}y},\,\, {\theta}_{xy}\big\{f(x,y)\}:=\frac{f(q^{-1}x,y)-f( x,qy)}{q^{-1}x-y},
\ee
which  turn out to be suitable for dealing with the Cauchy polynomials.  
%

Cao \cite{Jain2013} defined another homogeneous $q$-difference opeators
\be 
\label{qoperator}
\mathbb{T}(a,zD_{xy})=\sum_{k=0}^\infty \frac{(a;q)_k}{(q;q)_k} \left(z\,D_{xy}\right)^k,\quad\,\mathbb{E}(a, z\theta_{xy})=\sum_{k=0}^\infty \frac{(a;q)_k}{(q;q)_k}\, \left(-z\,\theta_{xy}\right)^k,
\ee
and obtain some results from the perspective of $q$-difference equations (see \cite{Jain2013}, for more details).  
 
 In order to reach our goals in this paper, we need the following Lemmas. 
 \begin{lemma}
  \label{thr}
  \bea
  \label{check1}
  \mathcal{C}_n^{(\alpha-n)}(x,y,b)&=&\mathbb{T} (q^{-\alpha}, bq^\alpha D_{xy}) \left\{p_n(x,y) \right\},\\
   \label{check2}
 \mathcal{D}_n^{(\alpha-n)}(x,y,b)&=&\mathbb{E} (q^{-\alpha} 
 ,bq^\alpha \theta_{xy})\left\{(-1)^nq^{-({}^n_2)} p_n(y,x) \right\}.
  \eea
 \end{lemma}
We now state  the
$q$-identities  asserted by Lemma \ref{lemma12} below.
 \begin{lemma}
\label{lemma12} It is asserted that
  \be 
\label{tO1}
\mathbb{T} (q^{-\alpha},  zD_{xy})\left\{\frac{(yt;q)_\infty}{(xt;q)_\infty}\right\}=\frac{(yt,q^{-\alpha} zt;q)_\infty}{(xt,zt;q)_\infty},\quad {\max\{|xt|,|zt|\} <1}
\ee 
and 
\be 
\label{tO2}
\mathbb{E}(q^{-\alpha}, z\theta_{xy})\left\{\frac{(xt;q)_\infty}{(yt;q)_\infty}\right\}=\frac{(xt,q^{-\alpha} zt;q)_\infty}{(yt,zt;q)_\infty},\quad  \max\{|yt|,|zt|\} <1. 
\ee 
 \end{lemma}
 \section{Proof of Theorems \ref{exgend} and   \ref{Tsxt}  }
\label{prel3}
In this section, we use the homogeneous  $q$-operators  defined in (\ref{qoperator}) to  prove    Theorems \ref{exgend} and \ref{Tsxt}. 
 
First, we give the identities (\ref{bella}) and (\ref{tbella})  below,   which will be used later in order to derive the   Rogers type formulas and extended Rogers type formula for the Cigler's polynomials  $\mathcal{C}_{n }^{(\alpha-n)}(x,y,b)$ and $\mathcal{D}_{n }^{(\alpha-n)}(x,y,b).$ 
  \begin{lemma}
  \label{OO1} It is asserted that
\be 
\label{bella}
\mathbb{T}(q^{-\alpha}, b  D_{xy}) \left\{\frac{p_n(x,y)\,(ys;q)_\infty}{(ys;q)_n(xs;q)_\infty}\right\} \\
=s^{-n}\frac{ (ys,q^{-\alpha}bs;q)_\infty}{ (xs,bs;q)_\infty} {}_{3}\Phi_2\left[\begin{array}{r}q^{-n},xs,bs;
 \\\\
ys,q^{-\alpha}bs;
 \end{array} 
q; q\right],
\ee
where $ {\max\{|xs|,|bs|\}<1};$
\be 
\label{tbella}
\mathbb{E} (q^{-\alpha},b  \theta_{xy}) \left\{\frac{p_n(y,x)\,(xs;q)_\infty}{(xs;q)_n(ys;q)_\infty}\right\}  
\\=s^{-n}\frac{ (xs,q^{-\alpha}bs;q)_\infty}{ (ys,bs;q)_\infty} {}_{3}\Phi_2\left[\begin{array}{r}q^{-n},ys,bs;
 \\\\
xs,q^{-\alpha}bs;
 \end{array} 
q; q\right],
\ee
where $ {\max\{|ys|,|bs|\}<1}$.
 \end{lemma}
 We are in position to prove  Theorems  \ref{exgend} and \ref{Tsxt}.
\begin{proof}[Proof  of Theorem  \ref{exgend}]  In light of (\ref{check1}), we have: 
\bea
 \sum_{n=0}^\infty\sum_{m=0}^\infty   \mathcal{C}_{n+m}^{(\alpha-n-m)}(x,y,b)\frac{  t^{n }}{(q;q)_n}\frac{  s^{m }}{(q;q)_m}&=&
\sum_{n=0}^\infty\sum_{m=0}^\infty  \mathbb{T}(q^{-\alpha},  q^{ \alpha}b  D_{xy}) \left\{ p_{n+m}(x,y)\right\}\frac{ t^{n }}{(q;q)_n}\frac{  s^{m }}{(q;q)_m}\cr
&= &\mathbb{T}(q^{-\alpha},  q^{ \alpha}b  D_{xy})\left\{\sum_{n=0}^\infty  p_{n}(x,y) \frac{  t^{n }}{(q;q)_n}\sum_{m=0}^\infty  p_{m }(x,yq^n)\frac{s^{m }}{(q;q)_m}\right\}
\cr
&= &\mathbb{T}(q^{-\alpha},  q^{ \alpha}b  D_{xy})\left\{\sum_{n=0}^\infty  \frac{  t^{n }}{(q;q)_n} \frac{p_{ n}(x,y)}{(ys;q)_n}\frac{(ys;q)_\infty}{(xs;q)_\infty}   \right\}\cr
&= &\sum_{n=0}^\infty  \frac{  t^{n }}{(q;q)_n}  \mathbb{T}(q^{-\alpha},  q^{ \alpha}b  D_{xy})\left\{\frac{p_{ n}(x,y)}{(ys;q)_n}\frac{(ys;q)_\infty}{(xs;q)_\infty}   \right\}. \nonumber
\eea
By using (\ref{bella}), we obtain:
\bea
 \sum_{n=0}^\infty\sum_{m=0}^\infty   \mathcal{C}_{n+m}^{(\alpha-n-m)}(x,y,b)\frac{  t^{n }}{(q;q)_n}\frac{  s^{m }}{(q;q)_m}&=& \frac{ (ys,bs;q)_\infty}{ (xs,q^{\alpha} bs;q)_\infty}\sum_{n=0}^\infty  \frac{ ( t/s)^{n }}{(q;q)_n}  {}_{3}\Phi_2\left[\begin{array}{r}q^{-n},xs,q^{\alpha} bs;
 \\\\
ys, bs;
 \end{array} 
q; q\right]\cr
&= &  \frac{ (ys,bs;q)_\infty}{ (xs,q^{\alpha} bs;q)_\infty}\sum_{n=0}^\infty \frac{(t/s)^n}{(q;q)_n}  \sum_{k=0}^\infty\frac{ (q^{-n},xs,q^{\alpha} bs;q)_k\,q^k}{(ys, bs,q;q)_k}
\cr
&= &\frac{ (ys,bs;q)_\infty}{ (xs,q^{\alpha} bs;q)_\infty}  \sum_{k=0}^\infty\frac{ (xs,q^{\alpha} bs;q)_k\,q^k}{(ys, bs,q;q)_k}  \sum_{n=0}^\infty \frac{ (q^{-n};q)_k(t/s)^n}{(q;q)_n} 
\cr
&= &\frac{ (ys,bs;q)_\infty}{ (xs,q^{\alpha} bs;q)_\infty}  \sum_{k=0}^\infty\frac{ (xs,q^{\alpha} bs;q)_k\,q^k}{(ys, bs,q;q)_k} 
\sum_{n=k}^\infty \frac{(t/s)^n (-1)^kq^{({}^k_2)-nk}}{(q;q)_{n-k}}
\cr
&= & \frac{ (ys,bs;q)_\infty}{ (xs,q^{\alpha} bs;q)_\infty}  \sum_{k=0}^\infty\frac{ (xs,q^{\alpha} bs;q)_k (-t/s)^k q^{-({}^k_2)}}{(ys, bs,q;q)_k}  \sum_{n=0}^\infty \frac{(tq^{-k}/s)^n}{(q;q)_{n}} \cr
&= &\frac{ (ys,bs;q)_\infty}{ (t/s,xs,q^{\alpha} bs;q)_\infty}  \sum_{k=0}^\infty\frac{ (xs,q^{\alpha} bs;q)_k}{(ys, bs,q;q)_k} \frac{(-t/s)^k q^{-({}^k_2)}}{(tq^{-k}/s;q)_k}\, \mbox{ by } (\ref{usu}) \cr
&= & \frac{ (ys,bs;q)_\infty}{ (t/s,xs,q^{\alpha} bs;q)_\infty}  \sum_{k=0}^\infty\frac{ (xs,q^{\alpha} bs;q)_k\,q^k}{(qs/t,ys, bs,q;q)_k}\cr
&= & \frac{ (ys,bs;q)_\infty}{ (t/s,xs,q^{\alpha} bs;q)_\infty}    {}_{4}\Phi_3\left[\begin{array}{r}xs,q^{\alpha} bs,0,0;
 \\\\
qs/t,ys, bs;
 \end{array} 
q; q\right].\nonumber
\eea
The proof of the   assertion (\ref{2exgen}) of Theorem \ref{exgend} is the same  to that of the first assertion (\ref{exgen}) by using the   representation  (\ref{check2}). 
The details involved are, therefore, being omitted here.
\end{proof}

\begin{proof}[Proof of Theorem \ref{Tsxt}] In view of
the formula (\ref{check1}),   we  observe that 
\bea
\label{BEL}
&&\sum_{n=0}^\infty\sum_{m=0}^\infty  \sum_{k=0}^\infty  \mathcal{C}_{n+m+k}^{(\alpha-n-m-k)}(x,y,b) \frac{  t^n\,s^m\omega^k}{(q;q)_{m}(q;q)_{n+m}(q;q)_k}\cr
& &\qquad=
\sum_{n=0}^\infty\sum_{m=0}^\infty  \sum_{k=0}^\infty\mathbb{T}(q^{-\alpha},  q^{ \alpha}b  D_{xy}) \left\{ p_{n+m+k}(x,y)\right\}\ \frac{  t^n\,s^m\omega^k}{(q;q)_{m}(q;q)_{n+m}(q;q)_k}
\cr
& &\qquad=\mathbb{T}(q^{-\alpha},  q^{ \alpha}b  D_{xy})\left\{\sum_{n=0}^\infty\sum_{m=0}^\infty\frac{  p_{n+m}(x,y) \,t^{n }  s^{m}}{(q;q)_m(q;q)_{n+m}} \sum_{k=0}^\infty  p_k(x,q^{n+m}y)\frac{ \omega^k}{(q;q)_k}   \right\}
\cr
& &\qquad=\mathbb{T}(q^{-\alpha},  q^{ \alpha}b  D_{xy})\left\{\sum_{n=0}^\infty\sum_{m=0}^\infty\frac{ p_{n+m}(x,y) t^{n }  s^{m}}{(q;q)_m(q;q)_{n+m}} \,\frac{(y\omega q^{n+m};q)_\infty}{(x\omega;q)_\infty} \right\}
\cr
& &\qquad= \sum_{n=0}^\infty  \sum_{m=0}^\infty \frac{t^{n }  s^{m}}{(q;q)_m(q;q)_{n+m}}    \mathbb{T}(q^{-\alpha},  q^{ \alpha}b  D_{xy})\left\{\frac{p_{n+m}(x,y) }{(y\omega ;q)_{n+m} }\frac{ (y\omega ;q)_\infty}{ (x\omega;q)_\infty}   \right\}. 
\eea
Setting $n$ by $n+m$ in  (\ref{bella}), we then obtain 
\be 
  \mathbb{T} (q^{-\alpha}, q^{ \alpha}b  D_{xy}) \left\{\frac{p_{n+m}(x,y) }{(y\omega ;q)_{n+m} }\frac{ (y\omega ;q)_\infty}{ (x\omega;q)_\infty}   \right\} 
=\omega^{-n-m}\frac{ (y\omega,b\omega;q)_\infty}{ (x\omega,q^{ \alpha} b\omega;q)_\infty} {}_{3}\Phi_2\left[\begin{array}{r}q^{-n-m},x\omega,q^{ \alpha} b\omega;
 \\\\
y\omega, b\omega;
 \end{array} 
q; q\right] 
\ee
which, in conjunction with (\ref{BEL}), gives
\bea
&&\sum_{n=0}^\infty\sum_{m=0}^\infty  \sum_{k=0}^\infty  \mathcal{C}_{n+m+k}^{(\alpha-n-m-k)}(x,y,b) \frac{  t^n\,s^m\omega^k}{(q;q)_{m}(q;q)_{n+m}(q;q)_k}\cr
&&\qquad=\frac{ (y\omega,b\omega;q)_\infty}{ (x\omega,q^{ \alpha} b\omega;q)_\infty} \sum_{n=0}^\infty  \sum_{m=0}^\infty \frac{(t/\omega)^{n }( s/\omega)^{m}}{(q;q)_m(q;q)_{n+m}}   {}_{3}\Phi_2\left[\begin{array}{r}q^{-n-m},x\omega,q^{ \alpha} b\omega;
 \\\\
y\omega, b\omega;
 \end{array} 
q; q\right] \cr
&&\qquad=\frac{ (y\omega,b\omega;q)_\infty}{ (x\omega,q^{ \alpha} b\omega;q)_\infty} \sum_{n=0}^\infty  \sum_{m=0}^\infty \frac{(t/\omega)^{n }( s/\omega)^{m}}{(q;q)_m(q;q)_{n+m}} \sum_{k=0}^{n+m}\frac{(q^{-n-m},x\omega,q^{ \alpha} b\omega;q)_k\,q^k}{(q,y\omega, b\omega;q)_k}  \cr
&&\qquad=\frac{ (y\omega,b\omega;q)_\infty}{ (x\omega,q^{ \alpha} b\omega;q)_\infty}\sum_{n=0}^\infty  \sum_{m=0}^\infty \sum_{n+m=k}^{\infty} \frac{(t/\omega)^{n }(s/\omega)^m (-1)^kq^{({}^k_2)-k(m+n)}}{(q;q)_m(q;q)_{n+m-k}} \frac{(x\omega,q^{ \alpha} b\omega;q)_k\,q^k}{(q,y\omega, b\omega;q)_k} \cr
&&\qquad=\frac{ (y\omega,b\omega;q)_\infty}{ (x\omega,q^{ \alpha} b\omega;q)_\infty} \sum_{m=0}^\infty\sum_{j=0}^\infty  \sum_{k=0}^\infty\frac{(s/t)^m}{(q;q)_m}\frac{(t/\omega)^{j+k}}{(q;q)_{j}}  \frac{(-1)^kq^{({}^k_2)-k(k+j)}  }{ (q;q)_k}
  \frac{(x\omega,q^{ \alpha} b\omega;q)_k\,q^k}{(y\omega, b\omega;q)_k} 
\cr
&&\qquad=\frac{ (y\omega,b\omega;q)_\infty}{ (x\omega,q^{ \alpha} b\omega;q)_\infty} \sum_{m=0}^\infty \frac{(s/t)^m}{(q;q)_m} \sum_{k=0}^\infty\frac{ (x\omega,q^{ \alpha} z\omega;q)_k(-t/\omega)^{k}q^{-({}^k_2) }}{(y\omega, b\omega,q;q)_k} \sum_{j=0}^\infty \frac{(tq^{-k}/\omega)^{j }}{(q;q)_{j}}   
\cr
&&\qquad=\frac{ (y\omega,b\omega;q)_\infty}{ (s/t,x\omega,q^{ \alpha} b\omega;q)_\infty}   \sum_{k=0}^\infty\frac{ (x\omega,q^{ \alpha} b\omega;q)_k}{(y\omega, b\omega,q;q)_k}    \frac{ (-t/\omega)^{k} q^{-({}^k_2) } }{(tq^{-k}/\omega;q)_\infty} \, \mbox{ by } (\ref{usu}) 
\cr
&&\qquad=\frac{ (y\omega,b\omega;q)_\infty}{ (s/t,t/\omega,x\omega,q^{ \alpha} b\omega;q)_\infty} 
  \sum_{k=0}^\infty\frac{  (x\omega,q^{ \alpha} b\omega;q)_k\, q^k }{(y\omega, b\omega,q\omega/t,q;q)_k}   \cr
&&\qquad=\frac{ (y\omega,b\omega;q)_\infty}{ (s/t,t/\omega,x\omega,q^{ \alpha} b\omega;q)_\infty} 
     {}_{4}\Phi_3\left[\begin{array}{r}x\omega,q^{ \alpha} b\omega,0,0; 
 \\\\
y\omega, b\omega,q\omega/t;
 \end{array} 
q; q\right].  \nonumber
\eea
The proof of the   assertion (\ref{2gextend}) of Theorem \ref{Tsxt} is the same  to that of the first assertion (\ref{gextend}) by using the   representation  (\ref{check2}). 
The details involved are, therefore, being omitted here.
\end{proof}
\section{Proof of Theorem \ref{TA1}}
\label{qdiffff}

The Hahn polynomials \cite{Hahn049,Hahn49}
 (or Al-Salam and Carlitz polynomials \cite{AlSalam,Cao2012A}) are defined as 
\be 
\label{ALSALAM}
\phi_n^{(a)}(x|q)=\sum_{k=0}^n \qbinomial{n}{k}{q} (a;q)_k\, x^k,\quad \psi_n^{(a)}(x|q)=\sum_{k=0}^n \qbinomial{n}{k}{q}q^{k(k-n)} (aq^{1-k};q)_k\,x^k.
\ee 
Srivastava and Agarwal \cite{SrivastavaAgarwal} utilized the method of transformation theory to deduce the following results, while Cao   \cite{Cao2012A} used the technique of exponential operator decomposition. For more information, please refer to \cite{Hahn049,Hahn49,AlSalam,SrivastavaAgarwal,Cao2009A,Cao2010A,Cao2012A}.
  \begin{lemma}(\cite[Eq. (3.20)]{SrivastavaAgarwal} and \cite[Eq. (5.4)]{Cao2012A}) 
\label{LEMMA41}
\be 
\label{21sums}
\sum_{n=0}^\infty \phi_n^{(\alpha)}(x|q) (\lambda;q)_n\frac{  t^n}{(q;q)_n} \\
= \frac{(\lambda t; q)_\infty }{(t;q)_\infty}   {}_2\Phi_1\left[
\begin{array}{rr} \lambda, \alpha;\\\\
 \lambda  t; \end{array}\,q; xt   
\right], \;{\max\{|t|, |xt|\}<1},
\ee
\be 
\label{c1sums}
\sum_{n=0}^\infty \psi_n^{(\alpha)}(x|q) (1/\lambda;q)_n\frac{  (\lambda tq)^n}{(q;q)_n} \\
=  \frac{(xtq; q)_\infty }{(\lambda xtq;q)_\infty}   {}_2\Phi_1\left[
\begin{array}{rr}1/ \lambda,1/ (\alpha x);\\\\
 1/(\lambda  xt); \end{array}\,q; \alpha q  
\right],
\ee
where $ {\max\{|\lambda x tq|, |\alpha q|\}<1}.$
 \end{lemma}  
  
 In the proof of Theorem \ref{TA1},  the following $q$-Chu-Vandermonde  formula will be needed. 
 \begin{lemma}($q$-Chu-Vandermonde   \cite[  Eq. (II.7)]{GasparRahman}) 
 \be
 \label{male}
   {}_2\Phi_1\left[
\begin{array}{rr}q^{-n},a;\\\\
 c; \end{array}\,q; \frac{cq^n}{a}
\right] =\frac{(c/a;q)_n}{(c;q)_n}.
\ee
 \end{lemma}
  \begin{proof}[Proof of   Theorem \ref{TA1}] We observe that
 \bea
  \sum_{n=0}^\infty \phi_n^{(\alpha)}(x|q) \mathcal{C}_{n}^{(\beta-n )}(u,v,b)  \frac{ \,t^n}{(q;q)_n}  & =&
\sum_{n=0}^\infty   \phi_n^{(\alpha)}(x|q)   \mathbb{T} (q^{-\beta},  q^{ \beta}b  D_{uv}) \left\{ p_{n }(u,v) \right\}  \frac{  t^n }{(q;q)_n}
\cr
&= &\mathbb{T} (q^{-\beta},  q^{ \beta}b  D_{uv}) \left\{\sum_{n=0}^\infty \phi_n^{(\alpha)}(x|q) p_{n }(u,v) \frac{   t^{n }   }{(q;q)_n}     \right\}
\cr
&= &\mathbb{T} (q^{-\beta},  q^{ \beta}b  D_{uv}) \left\{\sum_{n=0}^\infty\sum_{n=k}^\infty  p_{n }(u,v) \frac{  (\alpha;q)_k\,x^k t^{n }   }{(q;q)_k(q;q)_{n-k}}     \right\}
\cr
&= &\mathbb{T} (q^{-\beta},  q^{ \beta}b  D_{uv}) \left\{\sum_{k=0}^\infty p_{k }(u,v) \frac{  (\alpha;q)_k\,(xt)^{k }   }{(q;q)_k }   \sum_{n=0}^\infty  p_{n }(u,v q^k) \frac{  t^{n }   }{(q;q)_k }     \right\}
\cr
&= &\mathbb{T} (q^{-\beta},  q^{ \beta}b  D_{uv}) \left\{\sum_{k=0}^\infty p_{k }(u,v) \frac{  (\alpha;q)_k\,(xt)^{k }   }{(q;q)_k }  \frac{(v  tq^k; q)_\infty }{(u t;q)_\infty}    \right\}
\cr
&= &
  \sum_{k=0}^\infty \frac{(\alpha;q)_k\, (x t)^{k } }{ (q;q)_{k }}   \mathbb{T} (q^{-\beta},  q^{ \beta}b  D_{uv})  \left\{\frac{p_k(u,v)}{(v t ;q)_k}\frac{(v t ;q)_\infty}{(u t;q)_\infty} \right\}. \nonumber
\eea
By using (\ref{bella}), we have:
\bea
 \sum_{n=0}^\infty \phi_n^{(\alpha)}(x|q) \mathcal{C}_{n}^{(\beta-n )}(u,v,b)  \frac{ \,t^n}{(q;q)_n} & = &
  \sum_{n=0}^\infty \frac{(\alpha;q)_n\, (x t)^{n } }{ (q;q)_{n }}  t^{-n}\frac{ (v t,  bt;q)_\infty}{ (u t, q^{ \beta} bt;q)_\infty} {}_{3}\Phi_2\left[\begin{array}{r}q^{-n},u t, q^{ \beta} bt;
 \\\\
v t, bt;
 \end{array} 
q; q\right]
\cr
&= &\frac{ (v t,  bt;q)_\infty}{ (u t, q^{ \beta} bt;q)_\infty} \sum_{n=0}^\infty \frac{(\alpha;q)_n\, x ^{n } }{ (q;q)_{n }}  {}_{3}\Phi_2\left[\begin{array}{r}q^{-n},u t, q^{ \beta} bt;
 \\\\
v t, bt;
 \end{array} 
q; q\right]
\cr
&= &\frac{ (v t,  bt;q)_\infty}{ (u t, q^{ \beta} bt;q)_\infty} \sum_{k=0}^n\frac{ (u t,q^{\beta} bt;q)_k\,q^k}{(v t, bt,q;q)_k}  \sum_{n=0}^\infty \frac{ (q^{-n};q)_k(\alpha;q)_n\, x ^{n }}{(q;q)_n} 
\cr
&= &\frac{ (v t,  bt;q)_\infty}{ (u t, q^{ \beta} bt;q)_\infty}  \sum_{k=0}^n\frac{ (u t,q^{\beta} bt;q)_k\,q^k}{(v t, bt,q;q)_k} 
\sum_{n=k}^\infty \frac{ (-1)^kq^{({}^k_2)-nk}\, (\alpha;q)_n\, x^n}{(q;q)_{n-k}}
\cr
&= & \frac{ (v t,  bt;q)_\infty}{ (u t, q^{ \beta} bt;q)_\infty}  \sum_{k=0}^n\frac{ (\alpha,u t,q^{\beta} bt;q)_k   (-x)^k q^{-({}^k_2) }}{(v t, bt,q;q)_k} \sum_{n=0}^\infty \frac{(\alpha q^k;q)_n\, (xq^{-k})^n}{(q;q)_{n}} \cr&= &\frac{ (v t,  bt;q)_\infty}{ (u t, q^{ \beta} bt;q)_\infty}  \sum_{k=0}^n\frac{ (\alpha,u t,q^{\beta} bt;q)_k }{(v t, bt,q;q)_k} (- x)^k q^{-({}^k_2) }\frac{(\alpha x;q)_\infty}{( q^{-k}x;q)_\infty}\cr
&= & \frac{ (\alpha x,v t,bt;q)_\infty}{ (x,u t,q^{\beta} bt;q)_\infty}  \sum_{k=0}^n\frac{ (\alpha,u t,q^{\beta} bt;q)_k }{(v t, bt,q;q)_k} \frac{(-x)^k q^{-({}^k_2) }}{(x q^{-k};q)_k}\,  \mbox{ by } (\ref{usu})\cr
&= & \frac{ (\alpha x,v t,bs;q)_\infty}{ (x,u t,q^{\beta} bt;q)_\infty}  \sum_{k=0}^n\frac{ (\alpha,\mu t,q^{\beta} bt;q)_k\,q^k}{(q/x,v t, bt,q;q)_k}\cr
&= & \frac{ (\alpha x,v t,bt;q)_\infty}{ (x,u t,q^{\beta} bt;q)_\infty}    {}_{4}\Phi_3\left[\begin{array}{r}\alpha,u t,q^{\beta} bt,0;
 \\\\
q/x,v t, bt;
 \end{array} 
q; q\right].\nonumber
\eea
The proof of the assertion (\ref{1sums}) is thus completed. \\
Letting $(\lambda,t)=(v/u,tu)$ in equation (\ref{c1sums}),  we then obtain
\bea
\label{llms}
&&\sum_{n=0}^\infty \psi_n^{(\alpha)}(x|q)p_n(v,u) \frac{  (qt)^n}{(q;q)_n}\cr
&&\qquad =\frac{(uxtq; q)_\infty }{(v  xtq;q)_\infty}   {}_2\Phi_1\left[
\begin{array}{rr}1/ (\alpha x), u/v;\\\\
 1/(v  xt); \end{array}\,q; \alpha q  
\right] 
\cr
&&\qquad = \frac{ (uxtq; q)_\infty }{ (v  xtq;q)_\infty} \sum_{k=0}^\infty\frac{(1/(\alpha x);q)_k (\alpha q)^k}{(q;q)_k}\frac{(u/v;q)_k   }{ (1/(v  xt);q)_k }\mbox{ by } (\ref{male})\cr
&&\qquad =\frac{ (uxtq; q)_\infty }{ (v  xtq;q)_\infty} \sum_{k=0}^\infty\frac{(1/(\alpha x);q)_k (\alpha q)^k}{(q;q)_k}  {}_2\Phi_1\left[
\begin{array}{rr}q^{-k}, 1/ (uxt);\\\\
 1/(v  xt); \end{array}\,q; \frac{uq^k}{v} 
\right] \cr
&&\qquad = \sum_{k=0}^\infty\frac{(1/(\alpha x);q)_k (\alpha q)^k}{(q;q)_k}\sum_{n=0}^\infty\frac{(q^{-k};q)_n  \,q^{nk} }{(q;q)_n}
 \frac{(1/(uxt); q)_n    (uxtq; q)_\infty }{ (1/(v  xt);q)_n(v  xtq;q)_\infty}\left(\frac{u}{v}\right)^n \cr
&&\qquad = \sum_{k=0}^\infty\frac{(1/(\alpha x);q)_k (\alpha q)^k}{(q;q)_k}\sum_{n=0}^\infty\frac{(q^{-k};q)_n \, q^{nk}   }{(q;q)_n}
 \frac{   (uxtq^{1-n}; q)_\infty }{ (v  xtq^{1-n};q)_\infty}. 
\eea
We observe that 
 \bea
 &&\sum_{n=0}^\infty \psi_n^{(\alpha)}(x|q) \mathcal{D}_{n}^{(\beta-n )}(u,v,b)  \frac{(-1)^nq^{({}^{n+1}_{\,\,\,\,\,2})} \,t^n}{(q;q)_n}  \cr
&&\qquad =
\sum_{n=0}^\infty   \psi_n^{(\alpha)}(x|q)   \mathbb{E} (q^{-\beta},  q^{ \beta}b  \theta_{uv}) \left\{ p_{n }(v,u) \right\}  \frac{    (qt)^{n }  }{(q;q)_n}
\cr
&&\qquad =\mathbb{E} (q^{-\beta},  q^{ \beta}b  \theta_{uv}) \left\{\sum_{n=0}^\infty \psi_n^{(\alpha)}(x|q) p_{n }(v,u) \frac{   (qt)^{n }   }{(q;q)_n}     \right\} \mbox{ by } (\ref{llms})
\cr
&&\qquad =\mathbb{E} (q^{-\beta},  q^{ \beta}b  \theta_{uv}) \left\{\sum_{k=0}^\infty\frac{(1/(\alpha x);q)_k (\alpha q)^k}{(q;q)_k}\sum_{n=0}^\infty\frac{(q^{-k};q)_n \, q^{nk}   }{(q;q)_n}
 \frac{   (uxtq^{1-n}; q)_\infty }{ (v  xtq^{1-n};q)_\infty}  \right\}
\cr
&&\qquad =\sum_{k=0}^\infty\frac{(1/(\alpha x);q)_k (\alpha q)^k}{(q;q)_k}\sum_{n=0}^\infty\frac{(q^{-k};q)_n \, q^{nk}   }{(q;q)_n}\mathbb{E} (q^{-\beta},  q^{ \beta}b  \theta_{uv}) \left\{  \frac{   (uxtq^{1-n}; q)_\infty }{ (v  xtq^{1-n};q)_\infty}    \right\}
\cr
&&\qquad =\sum_{k=0}^\infty\frac{(1/(\alpha x);q)_k (\alpha q)^k}{(q;q)_k}\sum_{n=0}^\infty\frac{(q^{-k};q)_n \, q^{nk}   }{(q;q)_n}   \frac{   (uxtq^{1-n},  bxtq^{1-n}; q)_\infty }{ (v  xtq^{1-n},  bxtq^{\beta+1-n};q)_\infty} 
\cr
&&\qquad = \frac{   (uxtq,  bxtq; q)_\infty }{ (v  xtq,  bxtq^{1+\beta};q)_\infty} \sum_{k=0}^\infty\frac{(1/(\alpha x);q)_k (\alpha q)^k}{(q;q)_k}\cr
&&\qquad\times \sum_{n=0}^\infty   \frac{ (q^{-k},uxtq^{1-n},  bxtq^{1-n};q)_n \, q^{nk} }{   (v  xtq^{1-n},  bxtq^{\beta+1-n},q; q)_n }
  \cr
&&\qquad = \frac{   (uxtq,  bxtq; q)_\infty }{ (v  xtq,  bxtq^{1+\beta};q)_\infty}\sum_{n=0}^\infty   \frac{(-1)^nq^{({}^n_2)} (1/(uxt),  1/(bxt);q)_n   }{   (1/(v  xt),  q^{-\beta }/(bxt),q; q)_n }\left(\frac{u q^{-\beta}}{v}\right)^n\cr
&&\qquad \sum_{k=n}^\infty\frac{(1/(\alpha x);q)_k (\alpha q)^k}{(q;q)_{k-n}}
 \cr
&&\qquad = \frac{   (uxtq,  bxtq; q)_\infty }{ (v  xtq,  bxtq^{1+\beta};q)_\infty}\sum_{n=0}^\infty   \frac{(-1)^nq^{({}^n_2)} (1/(\alpha x),1/(uxt),  1/(bxt);q)_n   }{   (1/(v  xt),  q^{-\beta }/(bxt),q; q)_n }\cr
&& \qquad   \times \left(\frac{\alpha u q^{1-\beta}}{v}\right)^n \sum_{k=0}^\infty\frac{(q^n/(\alpha x);q)_k (\alpha q)^k}{(q;q)_{k}}\cr
&&\qquad = \frac{   (uxtq,  bxtq; q)_\infty }{ (v  xtq,  bxtq^{1+\beta};q)_\infty}\cr
&&\qquad\times \sum_{n=0}^\infty   \frac{(-1)^nq^{({}^n_2)} (1/(\alpha x),1/(uxt),  1/(bxt);q)_n   }{   (1/(v  xt),  q^{-\beta }/(bxt),q; q)_n }\left(\frac{\alpha u q^{1-\beta}}{v}\right)^n \frac{(q^{1+n}/  x;q)_\infty }{(\alpha q;q)_{\infty}}
\cr
&&\qquad = \frac{   (q/x,uxtq,  bxtq; q)_\infty }{ (\alpha q,v  xtq,  bxtq^{1+\beta};q)_\infty}\cr
&&\qquad \times \sum_{n=0}^\infty   \frac{(-1)^nq^{({}^n_2)} (1/(\alpha x),1/(uxt),  1/(bxt);q)_n   }{   (q/x,1/(v  xt),  q^{-\beta }/(bxt),q; q)_n }\left(\frac{\alpha u q^{1-\beta}}{v}\right)^n \cr
&&\qquad = \frac{   (q/x,uxtq,  bxtq; q)_\infty }{ (\alpha q,v  xtq,  bxtq^{1+\beta};q)_\infty}  {}_{3}\Phi_3\left[\begin{array}{rr}1/(\alpha x),1/(uxt),  1/(bxt);
 \\\\
q/x,1/(v  xt),  q^{-\beta }/(bxt);
 \end{array} 
q; \frac{\alpha u q^{1-\beta }}{v} \right], \nonumber
\eea
 which achieves the proof of the assertion (\ref{2sums}) of Theorem \ref{TA1}.  
 \end{proof}
\section{Another Srivastava-Agarwal type generating functions for the  Al-Salam-Carlitz polynomials}
\label{q5}
  In this section, we  derive another Srivastava-Agarwal type generating functions for the  Al-Salam-Carlitz polynomials (\ref{ALSALAM}).  As an application of Srivastava-Agarwal type generating functions,    we deduce two   interesting   transformation formulas between  ${}_2\Phi_1, \,  {}_2\Phi_2$  and ${}_3\Phi_2$. 
\begin{thm} 
\label{corld} For  $\alpha\in\mathbb{R}$, we have:
\be 
\label{1ss}
\sum_{n=0}^\infty  \phi_n^{(\alpha)}(x|q)  (\lambda;q)_n  \frac{ \,t^n}{(q;q)_n}  
= \frac{ (\alpha x,\lambda t;q)_\infty}{ (x,  t;q)_\infty}    {}_{3}\Phi_2\left[\begin{array}{r}\alpha,  t,0;
 \\\\
q/x,\lambda t;
 \end{array} 
q; q\right],\,  {\max\{|t|,|x|\} <1},
\ee
\be 
\label{2ss}
\sum_{n=0}^\infty  \psi_n^{(\alpha)}(x|q) (1/\lambda;q)_n  \frac{  (\lambda tq)^n}{(q;q)_n} 
= \frac{   (q/x, xtq; q)_\infty }{ (\alpha q, \lambda xtq ;q)_\infty}  {}_{2}\Phi_2\left[\begin{array}{rr}1/(\alpha x),1/(xt);\\
 \\
q/x,1/(\lambda xt);
 \end{array} 
q; \frac{\alpha   q}{\lambda } \right],
\ee
where $ {\max\{|\alpha q|,|xt|\} <1}$.
\end{thm}
We now derive two transformation formulas for $q$-series.
\begin{thm}
\label{ttL} We have:
\bea 
\label{1sdss}
   {}_2\Phi_1\left[
\begin{array}{rr}a, b;\\\\
c; \end{array}\,q; z   
\right] = \frac{ (abz/c;q)_\infty}{ (az/c;q)_\infty}    {}_{3}\Phi_2\left[\begin{array}{r}b,  c/a,0;
 \\\\
qc/(az),c;
 \end{array} 
q; q\right].
\eea
\end{thm}
\begin{thm}\cite[Eq. (III.4)]{GasparRahman}
\label{Coroy2L} We have:
\bea 
\label{2sdss}
   {}_2\Phi_1\left[
\begin{array}{rr}a, b;\\\\
c; \end{array}\,q; z   
\right]= \frac{   (bz; q)_\infty }{ (z;q)_\infty}  {}_{2}\Phi_2\left[\begin{array}{rr}b,c/a;\\
 \\
bz,c;
 \end{array} 
q; az\right].
\eea
\end{thm}
      \begin{remark}
   For $b=0$ and $u=1$,   the assertion (\ref{1sums}) of Theorem \ref{TA1}    reduces to    (\ref{1ss}). For $b=0,\,u=1/\lambda,\,v=1$ and $t=\lambda t$,    the assertion  (\ref{2sums}) of Theorem \ref{TA1}    reduces to    (\ref{2ss}).
    \end{remark}
\begin{remark}
  Comparing the assertion  (\ref{21sums}) of Lemma \ref{LEMMA41} and the assertion  (\ref{1ss}) of Theorem \ref{corld} and upon setting $\lambda=a,\,\alpha=b,\,\lambda t=c$ and $xt=z$, we obtain (\ref{1sdss}).
    \end{remark}
    \begin{remark}
 Comparing the assertion  (\ref{c1sums}) of Lemma \ref{LEMMA41} and the assertion  (\ref{2ss}) of Theorem \ref{corld} and upon setting $
  1/ \lambda=a,1/ (\alpha x)=b,\,  
 1/(\lambda  xt)=c$ and $  \alpha q=z$, we then obtain (\ref{2sdss}).
    \end{remark}

\section{Concluding remarks and observations}
\label{conclusion}
 
In our present investigation, we have used  two $q$-operators  $
\mathbb{T}(a,zD_{xy})$ and $\mathbb{E}(a, z\theta_{xy})$
to derive  Rogers formulas,    extended Rogers formulas and Srivastava-Agarwal   type bilinear generating functions  for Cigler's polynomials  by means of the $q$-difference equations.
We have also briefly described relevant 
connections of various special cases 
and consequences of our main results 
with several known results.
 
It is believed that the $q$-series   identities, which we
have presented in this paper, as well as the various related recent works
cited here, see also \cite{HMS_JC_SA2020}, will provide encouragement and motivation for further researches
on the topics that are dealt with and investigated in this paper. 

\medskip

\section*{Acknowledgments}
The author would like to thank the referees and editors for their many valuable comments and suggestions. This work was supported by the Zhejiang Provincial Natural Science Foundation of China (No.~LY21A010019).

\end{CJK*}
\end{document}